\documentclass[preprint,12pt]{elsarticle}

\usepackage{enumerate}
\usepackage{amsmath}
\usepackage{amsfonts}
\usepackage{amscd}
\usepackage{bm} 

\usepackage{ amsmath,amsfonts,amssymb,amscd,amsthm }
\usepackage{booktabs}
\usepackage{multirow}

\usepackage[usenames]{color}
\allowdisplaybreaks  

\usepackage{graphicx}
\usepackage{subfigure}
\usepackage{caption}
\newtheorem{thm}{Theorem}[section]

\newtheorem{pro}{Proposition}[section]
\newtheorem{lem}{Lemma}[section]
\newtheorem{cor}{Corollary}[section]

\newtheorem{assumption}{Assumption}

\renewcommand{\P}{\mathbb{P}}
\newcommand{\E}{\mathbb{E}}
\newcommand{\R}{\mathbb{R}}
\newcommand{\N}{\mathbb{N}}

\begin{document}

\begin{frontmatter}

\title{Asymptotic Properties of  a Branching-type Overloaded Polling Network}

\author[1]{Zaiming Liu}
\ead{math\_lzm@csu.edu.cn}
\author[1]{Yuejiao Wang}
\ead{wangyuejiaohujing@163.com}
\author[2]{Yuqing Chu\corref{cor1}}
\ead{chuyuqing@whut.edu.cn}
\author[3]{Yingqiu Li}
\ead{liyq-2001@163.com}

\cortext[cor1]{Corresponding author}
\address[1]{Department of Mathematics and Statistics, Central South University, Changsha 410083, Hunan, PR China}
\address[2]{School of Science, Wuhan University of Technology, Wuhan 430070, Hubei, PR China}
\address[3]{School of Mathematics and Statistics, Changsha University of Science and Technology, Changsha 410004, Hunan, PR China}

\begin{abstract}
 In this paper, we consider an $N$-queue overloaded polling network attended by a single cyclically roving server. Upon the completion of his service, a customer is either routed to another queue or leaves the system.  All the switches are instantaneous and random multi-gated service discipline is employed within each queue.
With the asymptotic theorem of multi-type branching processes and the exhaustiveness of service discipline, the fluid asymptotic process of the scaled joint queue length process is investigated. The fluid limit is of the similar shape with that of the polling system without rerouting policy, which allows us to optimize the gating indexes. Additionally, a stochastic simulation is undertaken to demonstrate the fluid limit and the optimization of the gating indexes to minimize the total population is considered.
\end{abstract}


\end{frontmatter}

\renewcommand{\thefootnote}{\fnsymbol{footnote}}

\section{Introduction}\label{intro_descri}
A typical polling system consists of a number of queues and a single server that visits the queues in a fixed order. In recent years, polling systems have become a fascinating area of research due to their wide applications in production-inventory system, air and railway transportation, the public health system, maintenance system, computer-communication system and flexible manufacturing system (see \cite{boon2011applications, levy1990polling} for overviews).
Along with that, a wide range of polling models have emerged.

In this paper, we consider a cyclic N-queue ($Q_1, \cdots, Q_N, N \geq 2$) polling system with random multi-gated service discipline within each queue and customer re-routing policies: after completing service at $Q_i$, a customer is either routed to $Q_j$ with probability $p_{i,j}$ or leaves the system with probability $p_{i,0}$. The possibility for re-routing of customers further enhances the already-extensive modeling capabilities of polling models, since in many applications, customers require service at more than one facility of the system. Actually, the models of customer re-routing arise naturally in various models of computer, communication and robotic systems (see for example \cite{boon2011queueing, boon2013waiting, sidi1990customer})
). One obvious example is a local area network in which terminals are interconnected in either a physical or logical structure (see \cite{sidi1992queueing}).

In the vast majority of papers that have appeared on polling models, it is almost invariably assumed that the system is stable and the stable performance measures are then concerned. With the advent of the era of Internet+, the study of critically or strictly super-critically loaded polling systems is vigorously pioneered due to the overloaded Internet channel or online shopping orders.

The heavy traffic ($\rho \to 1$, $\rho$ is the load of the system) behaviors  have gained an ascending attention in the last two decades pioneered by Coffman et al. \cite{coffman1995polling,coffman1998polling}.
By utilizing the connection with multi-type branching process,
van der Mei \cite{van2007towards} considered a unifying theory on branching-type polling models under heavy-traffic assumptions.
In the similar way, Boon et al. \cite{boon2011queueing} discussed the heavy-traffic asymptotic behaviors of a gated polling system with customer re-route policy.
Furthermore, Liu et al. \cite{liu2015asymptotic} extended the results in \cite{boon2011queueing} to the analogous system with a general branching-type service policy in the same form.
As for a non-branching type polling system, one example is that Liu et al. \cite{liu2016heavy}  investigated the heavy-traffic behavior of a priority polling system consisting of three M/M/1 queues with threshold policy and proved  that the scaled queue-length of the critically loaded queue is exponentially distributed, independent of that of the stable queues.

The study of overloaded ($\rho > 1$) service system is important to control or predict how fast it blows up over time. However, hardly any attention has been given to the overloaded polling system.  The few literature refers to \cite{puha2006fluid,talreja2008fluid,jennings2012overloaded,remerova2013random}.
By using measure-valued state descriptor, Puha et al. \cite{puha2006fluid} proved that the overloaded $GI/GI/1$ processor sharing queues converge in distribution to supercritical fluid models and a fluid limit result is proved as first order approximations to overloaded processor sharing queues.
Using both fluid and diffusion limits, Jennings et al. \cite{jennings2012overloaded} showed that the virtual waiting time process of an overloaded Multi-class FIFO(first-in-first-out) queue with abandonments converges to a limiting deterministic fluid process.
Instead Remerova et al. \cite{remerova2013random} showed the fluid asymptotic process for the joint queue length process on an overloaded branching-type cyclic polling system by using the asymptotic properties of multi-type branching processes.

In the present paper, we dedicate to the investigation of the overloaded asymptotic fluid process of the scaled joint queue length process. The fluid model associated with heavily loaded polling network is contained in \cite{liu2015asymptotic}.
The work carried out here is a natural progression from \cite{liu2015asymptotic} and a natural extension of \cite{remerova2013random}. 
due to the exhaustiveness, the fluid process here has the same shape with that in \cite{remerova2013random}.
For the linear shape and constant growth rate property of the fluid limit, we could optimize the system through minimizing the total queue length by choosing appropriate parameters (gating indexed and rerouting probabilities for example).


The rest of the paper is organized as follows. In Section \ref{main-result}, we describe precisely the polling model. In Sections \ref{multitype-branchin} and  \ref{connection-polling-branchin}, we introduce the super-critical property of multi-type branching process and construct the multi-branching process associated with our polling model respectively. Section \ref{fluid-main-results} provides the main results, where Lemmas \ref{eq-variable} and \ref{ineq-infty} are proved in Section \ref{prove-lemma} and Theorems \ref{thmQ} and \ref{thmcheckQ} are proved in Section \ref{prove-theorem}. Section \ref{numerical} discusses some numerical issues  including the stochastic simulation to test Theorems \ref{thmQ} and \ref{thmcheckQ} and the optimization of gating indexes to minimize the average growth rate of the total population. Section \ref{conclude}  concludes and provides an outlook on potential further research of our paper.

\section{Model description }\label{main-result}

  Consider an asymmetric cyclic polling model that consists of $N \geq 2$ queues, $Q_1, \cdots, Q_N$, and a single server that visits the queues in a cyclic order.
Customers arrive at $Q_i$ according to a Poisson process $E_i(\cdot)$ with rate $\lambda_i$. The service time of each customer at $Q_i$ is a random variable $B_i$ with finite mean value $\E B_i = 1/\mu_i$. The service discipline at $Q_i$ is multi-gated with gating index of a random variable $\kappa_i \in \N\cup \{\infty\}$ (denoted by $\kappa_i$-gated, see {\cite{remerova2013random}}).
The interarrival times, the service times and the gating indices (for different queues and for different visits) are assumed to be mutually independent.

Upon completion of service at $Q_i, i = 1, \cdots, N$, a customer is either routed to $Q_j, j = 1, \cdots, N$ with probability $p_{i,j}$ or leaves the system with probability $p_{i,0}$, where
\begin{eqnarray*}
\sum_{i=1}^{N} p_{i,0} > 0 \quad \text{and} \quad \sum_{j=0}^{N}p_{i,j}=1.
\end{eqnarray*}
We assume that all the switches of customers or servers between queues are instantaneous. Additionally, when the system becomes empty, the server travels a full cycle and subsequently stops right before $Q_1$ until a new  arrival occurs and then cycles along the queues to serve that customer.

The total arrival rate at $Q_i$ is denoted by $\gamma_i$, which is the unique solution of the following set of linear equation (\cite{boon2013waiting}):
\begin{eqnarray*}
\gamma_i = \lambda_i + \sum_{j=1}^{N} \gamma_j p_{j,i} \qquad i = 1, \cdots, N.
\end{eqnarray*}
The offered load to $Q_i$ equals to $\rho_i =  \gamma_i/\mu_i$ and the total load equals to $\rho = \sum_{i=1}^{N} \rho_i$. Furthermore, we need two more assumptions.
\begin{assumption}
For all $i = 1, \cdots, N$, $\rho_i < 1$, and $\rho > 1$.
\end{assumption}

\begin{assumption}
For all $i = 1, \cdots, N$, $\E B_i \log B_i < \infty$.
\end{assumption}

Throughout the paper, we focus our attention on the overloaded behavior of the queue length process ${\bf X}(\cdot) = (X_1, \cdots, X_N)(\cdot)$,
where $X_i(t)$ is the number of customers at $Q_i$ at time $t$.

For simplicity, we adopt the following notations.
\begin{enumerate}
  \item Let $\N = \{0, 1, \cdots\}$ and $\N^+ = \{1, 2, \cdots\}$.
  \item The vector with all coordinates equal to 0 is denoted by ${\bf 0}$, with all coordinates equal to 1 by ${\bf 1}$, and with coordinate $i$ equals to $1$ and the other coordinates equal to 0 by ${\bf e}_i$.
  \item For vectors ${\bf x} = (x_1, \cdots, x_N)$ and ${\bf y} = (y_1, \cdots, y_N)$, define the following operations:

$\bullet$ coordinate-wise product: ${\bf x} \times {\bf y} = (x_1y_1, \cdots, x_N y_N)$;

$\bullet$ power: if all $x_i > 0, 1 \leq i \leq N$, then ${\bf x^y} = \prod_{i=1}^{N} x_i^{y_i}$;

$\bullet$ binomial coefficient: if ${\bf y } \leq {\bf x}$, then $\bigg( \begin{array}{c}
                                                                     {\bf x} \\
                                                                     {\bf y}
                                                                   \end{array}
\bigg) = \prod_{i=1}^{N}\Big(\begin{array}{c}
                           x_i \\
                           y_i
                         \end{array}
\Big) = \prod_{i=1}^N x_i!/$ $y_i!(x_i - y_i)!$.
\item If random objects $X$ and $Y$ are equal in distribution, we write $X \overset{d}{=} Y$ and say that $X$ is a copy of $Y$.
\end{enumerate}

\section{Multi-type branching process}\label{multitype-branchin}
We consider a general multi-type branching process (MTBP) with $N$ particle types, denoted by ${\bf Z}_{n} = (Z_n^{(1)}, \cdots, Z_n^{(N)}), n \in \N$, where $Z_n^{(i)}$ is the number of type-$i$ particles in the $n$th generation for $i = 1, \cdots, N$, $n \in \N$.
\begin{itemize}
\item Define the immigration distribution by
\begin{eqnarray*}
G({\bf k}): = \P({\bf Z}_{n+1} = {\bf k}|{\bf Z}_{n}=0), \quad {\bf k} \in \N^{N}.
\end{eqnarray*}
  \item Let ${\bf M} = \{m_{i,j}\}_{i,j=1}^{N}$ be the mean offspring matrix, where $m_{i,j}$ is the expected number of type $j$ offspring of a single type $i$ particle in one generation.
  \item  Let the vectors ${\bf u} = (u_1, \cdots, u_N)$ and ${\bf v} = (v_1, \cdots, v_N)$ be the right and left eigenvectors corresponding to the maximal  real-valued, positive eigenvalue $\theta$ of ${\bf M}$, commonly referred to as the maximum eigenvalue (\cite{athreya1972branching}), normalized such that ${\bf v u^\top} = 1$.
  \item Let $ {\bf q} = (q_1, \cdots, q_N)$, where $q_i$ is the probability of eventual extinction of the process initiated with a single particle of type $i$, i.e.,
$$q_i=\P({\bf Z}_n = {\bf 0} \ \text{for some} \ n |{\bf Z}_0 = {\bf e}_i).$$
\end{itemize}



It follows that the auxiliary process $\{{\bf Z}^{(n)}\}_{n \in \N}$ has the following asymptotics.

\begin{pro} (\cite{athreya1972branching})
Given ${\bf Z}^{(0)} = {\bf e}_i,$
\begin{eqnarray*}
\frac{{\bf Z}^{(n)}}{\theta^n} \to \xi_i {\bf v}\qquad \mbox{almost surely (a.s.) } \qquad \mbox{as} \; n \to \infty,
\end{eqnarray*}
where the distribution of the random variable $\xi_i$ has a point mass  $q_i < 1$ at 0 and a continuous density function on $(0, \infty)$ with $\E \xi_i = u_i$.
\end{pro}

\section{Branching property of multi-gated polling system}\label{connection-polling-branchin}
To start with, define $t^{(n)}$ as the time point that the server reaches right before $Q_1$ for the $n$th time and $t_i^{(n)}$ as the time point that the  server reaches $Q_i$ for the $n$th time ($n\in \N^+, i=1,2,\ldots,N$), then
\begin{enumerate}
  \item $t^{(n)} \leq t_1^{(n)} \leq \cdots \leq t_{N + 1}^{(n)} =t^{(n+1)}$;
  \item If the system is empty at $t^{(n)}$, then the interval $[t^{(n)}, t_1^{(n)})$ is the period of waiting until the first arrival, otherwise $t^{(n)} = t_1^{(n)}$;
  \item The interval $[t_i^{(n)}, t_{i+1}^{(n)})$ is the visit time at $Q_i$ following $t^{(n)}$, with $t_i^{(n)} = t_{i+1}^{(n)}$ if $Q_i$ is empty.
\end{enumerate}
Typically, we assume $t^{(1)}=0$, which means the system is empty at $t=0$.

{\bf Branching Property}{\cite{resing1993polling}} If the server arrives at $Q_i$ to find $k_i$ customers there, then during the course of the server's visit, each of these $k_i$ customers will effectively be replaced in an i.i.d. manner by a random population having probability generating function (p.g.f) $h_i(z_1,z_2,\ldots,z_N)$, which can be any $N$-dimensional p.g.f..

By Resing \cite{resing1993polling}, branching property implies that the queue length sequence $\{{\bf X}(t^{(n)})\}_{n \in \N}$
forms a multi-type branching process with immigration in state ${\bf 0}$. Then the probability for the process $\{{\bf X}(t^{(n)})\}_{n \in \N}$ to return to ${\bf 0}$ is given by
\begin{eqnarray*}
q_G: = \sum_{{\bf k} \in \N^N} G({\bf k}) {\bf q^k}.
\end{eqnarray*}
Subsequently, we give some notations associated with the branching-type polling system.
\begin{itemize}
  \item Define ${\bf \check{L}}_{i}=(\check{L}_{i,1}, \cdots, \check{L}_{i,N})$ as the visit offspring of a customer at $Q_i$, which equals in distribution to ${\bf X}(t_{i+1}^{(n)})$ given that ${\bf X}(t_i^{(n)}) = {\bf e}_i$ (its distribution does not depend on $n$).
  \item Define ${\bf{L}}_{i}:=(L_{i,1}, \cdots, L_{i,N})$ as the  session offspring of a customer at $Q_i$, which equals in distribution to ${\bf X}(t^{(n+1)})$ given that ${\bf X}(t^{(n)}) = {\bf e}_i$ (its distribution does not depend on $n$).
  \item Denote the mean visit offspring of a customer at $Q_i$ and the mean session offspring by ${\bf \check{m}}_i = (\check{m}_{i,1}, \cdots, \check{m}_{i,N})=\E {\bf \check{L}}_{i} $ and ${\bf m}_i = (m_{i,1}, \cdots, m_{i,N})= \E{\bf L }_{i}$, respectively.
\item Define $T_i$ as the visit duration at $Q_i$ which equals in distribution to
$t_{i+1}^{(n)} - t_{i}^{(n)}$ given that $X_i(t_i^{(n)})=1$. 
\end{itemize}
 Then the immigration distribution is given by
\begin{eqnarray}\label{definition-G(k)}
G({\bf k})= \frac{\sum_{i=1}^{N} \lambda_i \P({\bf L}_i = {\bf k})}{\sum_{i=1}^{N} \lambda_i}, \quad {\bf k} \in \N^{N}.
\end{eqnarray}

Before proceeding further, we give a lemma on the exhaustiveness $f_i$ of the service descipline at $Q_i$, which is defined by (see \cite{van2000polling})
\begin{equation*}
  f_i = 1 - \frac{\partial}{\partial z_i}h_i(z_1,z_2, \cdots, z_N)|_{{\bf z} = {\bf 1}}=1 - \E \check{L}_{i,i}.
\end{equation*}
 Virtually, it has an appealing interpretation: during the course of the server's visit at $Q_i$, each customer present at the start of the visit to $Q_i$ will be replaced by a number of customers with mean $1-f_i$ at the end of the visit to $Q_i$.

\begin{lem}\label{exhaustiveness}
In our model, for $Q_i$, we have
\begin{align}
  &f_i=1-\E (\frac{\lambda_i}{\mu_i} + p_{i,i})^{\kappa_i},\label{exhuastive}\\
&t_i :=\E T_i=\frac{1 - \E (\frac{\lambda_i}{\mu_i} + p_{i,i})^{\kappa_i}}{\mu_i(1 - \frac{\lambda_i}{\mu_i}-p_{i,i})}.\label{visit}
\end{align}
\end{lem}


Let $B_i^E$ and $Y_i$ be the total service time and the number of services of a customer in $Q_i$ before he is either routed to $Q_j, j \neq i$ or leaves the system. Obviously, $\P(Y_i = n) = p_{i,i}^{n-1}(1-p_{i,i}), (n = 1, 2, \cdots)$ and
\begin{equation*}
b_i^E = \E B_i^E = \E \sum_{j=1}^{Y_i} B_{i,j} = \E Y_i \E B_i= \frac{1}{\mu_i(1 - p_{i,i})},
\end{equation*}
where $\{B_{i,j}, j=1, 2, \cdots\}$ are i.i.d. copies of $B_i$.
By Lemma 1 in \cite{liu2015asymptotic}, the mean offspring matrix ${\bf M}$ is  given in the following lemma.
\begin{lem}\label{eq-variable}
For the cyclic branching-type polling system, the mean matrix ${\bf M}$ is given by
\begin{equation*}\label{eq4.7}
    \mathbf{M}=\mathbf{M}_1\ldots \mathbf{M}_N,
  \end{equation*}
where $\mathbf{M}_k=\left(m^{(k)}_{i,j}\right)$ and
  \begin{equation*}
    m^{(k)}_{i,j}=\left\{
       \begin{array}{ll}
         \delta_{\{i=j\}}, & \hbox{$i\neq k$,} \\
         1-f_i, & \hbox{$i=k=j$,} \\
         f_i\varphi_i(\mu_ip_{i,j}+\lambda_j), & \hbox{$i=k\neq j$,}
       \end{array}
     \right.
  \end{equation*}
where $\delta_{F}$ denotes the indicator function on $F$ and $\varphi_i=\frac{b_i^E}{1-\lambda_ib_i^E}$.
 \end{lem}

Actually, $\mathbf{M}_k$ is the mean session offspring during the visit time on $Q_k$. Hence, for all $i$,
\begin{equation*}
   \check{m}_{i,j}=m^{(i)}_{i,j}=\left\{
       \begin{array}{ll}
         1-f_i, & \hbox{$i=j$,} \\
         f_i\varphi_i(\mu_ip_{i,j}+\lambda_j), & \hbox{$i\neq j$.}
       \end{array}
     \right.
  \end{equation*}
Thus, for the $m_{i,j}$, we also have the following recursive formula:
\begin{eqnarray*}
m_{N,j} = \check {m} _{N, j}, \qquad \mbox{ for all }\quad j,
\end{eqnarray*}
and, for $ i = 1, \cdots, N-1$, ${\bf m}_i$ is computed via ${\bf m}_{i+1}$,
\begin{eqnarray*}
m_{i,j} = \check{m}_{i,j}\delta_{\{i \geq j\}} + \sum_{k=i+1}^{N} \check{m}_{i,k}, \qquad \mbox{for all } \quad j.
\end{eqnarray*}

By the properties of the maximal eigenvalue $\theta$ of the mean matrix ${\bf M}$, 
the following Lemma \ref{extinc}  construct the connection between extinction probability and the maximal eigenvalue $\theta$ in  branching-type polling systems. For notational convenience, denote $\theta(\rho)$ by $\theta$.
\begin{lem}\label{extinc}
For the Perron-Frobenius eigenvalue $\theta$ and the extinction probabilities $q_i$, we have $\theta > 1$ and $q_i < 1$ for all $i$. By the latter, $q_G < 1$, too.
\end{lem}

It follows that the auxiliary process $\{{\bf X}(t^{(n)})\}_{n \in \N}$ has the following asymptotics.

\begin{pro}\label{pro-limit}
If the first arriving customer arrives at $Q_i$ after $t=0$, then
\begin{eqnarray*}
\frac{{\bf X}(t^{(n)})}{\theta^n} \to \xi_i {\bf v}\qquad \mbox{almost surely (a.s.) } \qquad \mbox{as} \; n \to \infty,
\end{eqnarray*}
where the distribution of the random variable $\xi_i$ has a jump of magnitude $q_i < 1$ at 0 and a continuous density function on $(0, \infty)$ and $\E \xi_i = u_i$.
\end{pro}

Lemma \ref{ineq-infty} below considers the finiteness of the corresponding moments for the offspring distribution of the multi-type branching processes $\{{\bf Q}(t^{(n)})\}_{n \in \N}$, which guarantees the non-degenerate of the random variable $\xi_i$.
\begin{lem}\label{ineq-infty}
For all $i$ and $j$, $\E L_{i,j} \log L_{i,j} < \infty$, where $0 \log 0: = 0$ by convention.
\end{lem}

\section{Fluid limit}\label{fluid-main-results}
To give the main results, two more notations are needed.
\begin{itemize}
  \item Let $\bar{B}_i$ be the total service time of a customer arriving at $Q_i$ from outside, $\bar{c}_i = \E \bar{B}_i$
 and  $A_{i}$ be the position of a customer  after completion of service at $Q_i$, for $i = 1, \cdots, N$, i.e.,
\begin{eqnarray*}\label{eq-def-Ai}
A_{i} = \left\{
            \begin{array}{ll}
              j, & \hbox{after receiving service at $Q_i$, a customer is routed to $Q_j$;} \\
              0, & \hbox{after receiving service at $Q_i$, a customer leaves the system.}
            \end{array}
          \right.
\end{eqnarray*}
Then $\P(A_i=j)=p_{i,j}$, $j =0, 1, \cdots, N$. By the law of total probability, we have
\begin{align*}
\bar{c}_i&=\E \bar{B}_i =  \E( \bar{B}_i|A_i =0) \P(A_i = 0) + \sum_{j=1}^{N} \E(\bar{B}_i | A_i = j) \P(A_i = j)\\
&= \E  B_i + \sum_{j=1}^{N} p_{i,j} \E \bar{B}_j= 1/\mu_i + \sum_{j=1}^{N} p_{i,j} \bar{c}_j.
\end{align*}
 It is also easy to deduce that $\rho = \sum_{i=1}^{N} \lambda_i \bar{c}_i$.
  \item For $n \in \N$, let
\begin{eqnarray*}
\eta_n: =\left\{
           \begin{array}{ll}
             \min\{k: t^{(k)} \geq \theta^n\}, & \hbox{ if $n \geq 0$;} \\
             0,  \;  &\;  \hbox{if $n < 0$.}
           \end{array}
         \right.
\end{eqnarray*}
\end{itemize}

\begin{thm}\label{thmQ}
There exist constants $\bar{b}_i \in(0, \infty)$ and ${\bf \bar{a}}_i = (\bar{a}_{i,1}, \cdots, \bar{a}_{i,N}) \in [0, \infty)^{N}, i=1, \cdots, N+1$, and a random variable $\xi$ with values in $[1, \theta)$ such that, for all $k \in \N$ and $i$,
\begin{eqnarray*}\label{eq-eta-q}
\frac{t_i^{(\eta_n + k)}}{\theta^n} \to \theta^k \bar{b}_i \xi  \qquad \mbox{and} \qquad \frac{X(t_i^{(\eta_n + k)})}{\theta^n} \to \xi \theta^k {\bf \bar{a}}_i \quad a.s. \quad \mbox{as} \quad n \to \infty.
\end{eqnarray*}
The $\bar{b}_i$ and ${\bf \bar{a}}_i$ are given by
\begin{equation*}\label{eq-def-ai-bi}
\bar{b}_1 = 1, \qquad \bar{b}_{i+1} = \bar{b}_i + \Big[\frac {v_i} {\alpha} + \lambda_i(\bar{b}_i - \bar{b}_1)+\sum_{j=1}^{i-1} p_{j,i} \mu_j(\bar{b}_{j+1} -\bar{ b}_j)\Big]t_i, \ \ \ i=1, \cdots, N;
\end{equation*}
and for $i = 1, \cdots, N$,
\begin{eqnarray*}\label{eq-def-bar-ai-equivalent}
{\bf \bar{a}}_1 = \frac{\bf v} {\alpha},
&&  \bar{a}_{i+1, j} = \left\{
                   \begin{array}{ll}
                    \bar{a}_{i,j} + [\lambda_j+ \mu_ip_{i,j}] (\bar{b}_{i+1} - \bar{b}_{i}) \qquad \qquad j \neq i, \\
                   \bar{a}_{i,i} + [\lambda_i -\mu_i(1-p_{i,i})] (\bar{b}_{i+1} - \bar{b}_i),  \ \quad j=i.
                   \end{array}
                 \right.
\end{eqnarray*}
where $\alpha =  \frac{\sum_{i=1}^N v_i \bar{c}_i}{\rho- 1}$.
For $x \in [1, \theta)$, the distribution of $\xi$ is given by
\begin{align*}
\P(\xi \geq x) =& \frac{1}{1- q_{G}}\sum_{\overset{{\bf k} \in \N^{N}} {|{\bf k}| \geq 1}} G({\bf k}) \sum_{\stackrel{{\bf l \leq k}}{|{\bf l}| \geq 1}} \big(\begin{array}{c}                                                                                                                                      {\bf k} \\                                                                                                                                     {\bf l}                                                                                                                                   \end{array}
\big){\bf(1- q)^l q^{k-l} } \\
&\times \P \bigg\{
                                                 \bigg\{\log_{\theta}
                                                             \big(\alpha \sum_{i=1}^{N}\sum_{j=1}^{l_i} \xi_i^{(j)}
                                                               \big)
                                                  \bigg \} \geq \log_{\theta} x
                                            \bigg\}
\end{align*}
where $\xi_i^{(j)}, j \in \N^+$ are i.i.d. copies of $\xi_i\delta_{\xi_i > 0}$.
\end{thm}

For each $n \in \N$, define the scaled queue length process
\begin{eqnarray}\label{def-sal-leng}
{\bf \bar{X}}^{(n)}(t) := \frac{{\bf X}(\theta^n t)}{\theta^n}, \qquad t \in [0, \infty).
\end{eqnarray}

\begin{thm}\label{thmcheckQ}
There exists a deterministic function ${\bf \bar{X}}(\cdot) = (\bar{X}_1, \cdots,\bar{X}_{N})(\cdot) \in [0, \infty)^N$
such that,
\begin{equation*}
{\bf \bar{X}}^{(n)}(\cdot) \to \xi {\bf \bar{X}}(\frac{\cdot}{\xi})  \quad a.s. \qquad \mbox{as} \quad n \to \infty,
\end{equation*}
uniformly on compact sets, where the random variable $\xi$ is defined in Theorem \ref{thmQ}.
For all $i = 1, \cdots, N$, the function ${\bf \bar{X}}(\cdot)$ is continuous and piecewise linear and specified by
\begin{eqnarray}\label{def-eq-barQ}
{\bf \bar{X}}(t) = \left\{
                     \begin{array}{ll}
                       0, & \hbox{if $t = 0$;} \\
                       \theta^k {\bf \bar{a}}_i + (t - \theta^k \bar{b}_i){\bm \lambda} +(t-\theta^k \bar{b}_i)\mu_i {\bf p}_i, & \hbox{if $t \in [\theta^k \bar{b}_i, \theta^k \bar{b}_{i+1})$},
                                                 \end{array}
                   \right.
\end{eqnarray}
where  ${\bm \lambda} = [\lambda_1, \cdots, \lambda_N]$ and ${\bf p}_{i}=[p_{i,1},\cdots,p_{i,i-1},p_{i,i}-1,p_{i,i+1},\cdots,p_{i,N}]$.
\end{thm}

\begin{figure}[!htb]
\centering
\includegraphics[width=4.5 in]{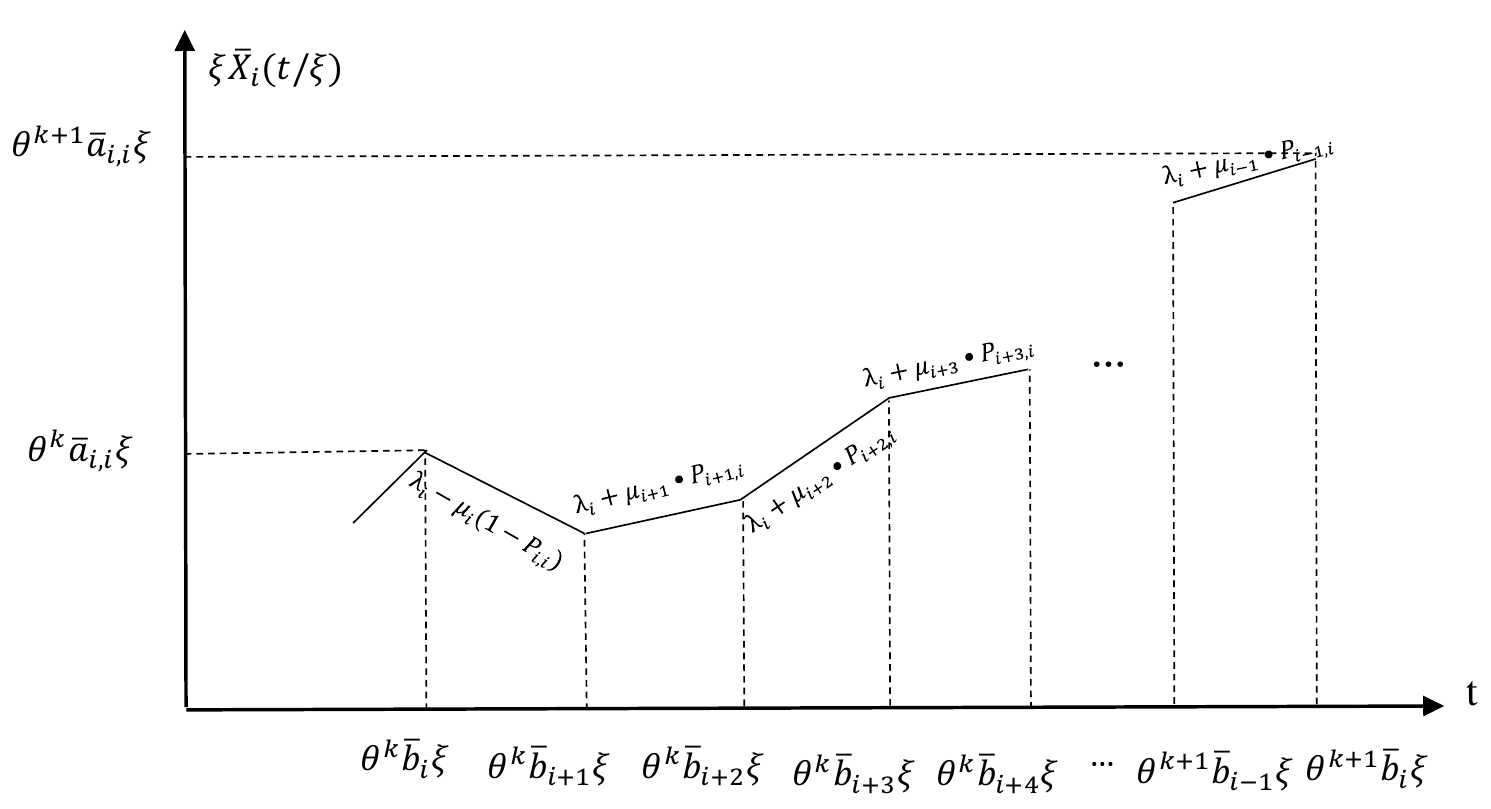}
\caption{Fluid limit of $Q_i$.}
\label{fig:graph}
\end{figure}

\begin{cor}\label{cor_total}
The limit total population $(\bar{X}_1 + \bar{X}_2 + \cdots + \bar{X}_N)(\cdot)$ grows at rate
$$(\lambda_1 + \cdots + \lambda_N) - p_{i,0} \mu_i$$
when $t \in [\theta^k \bar{b}_i, \theta^{k} \bar{b}_{i+1})$ for all $k \in \N$.
\end{cor}

\begin{figure}[!htb]
\centering
\includegraphics[width=3.8 in]{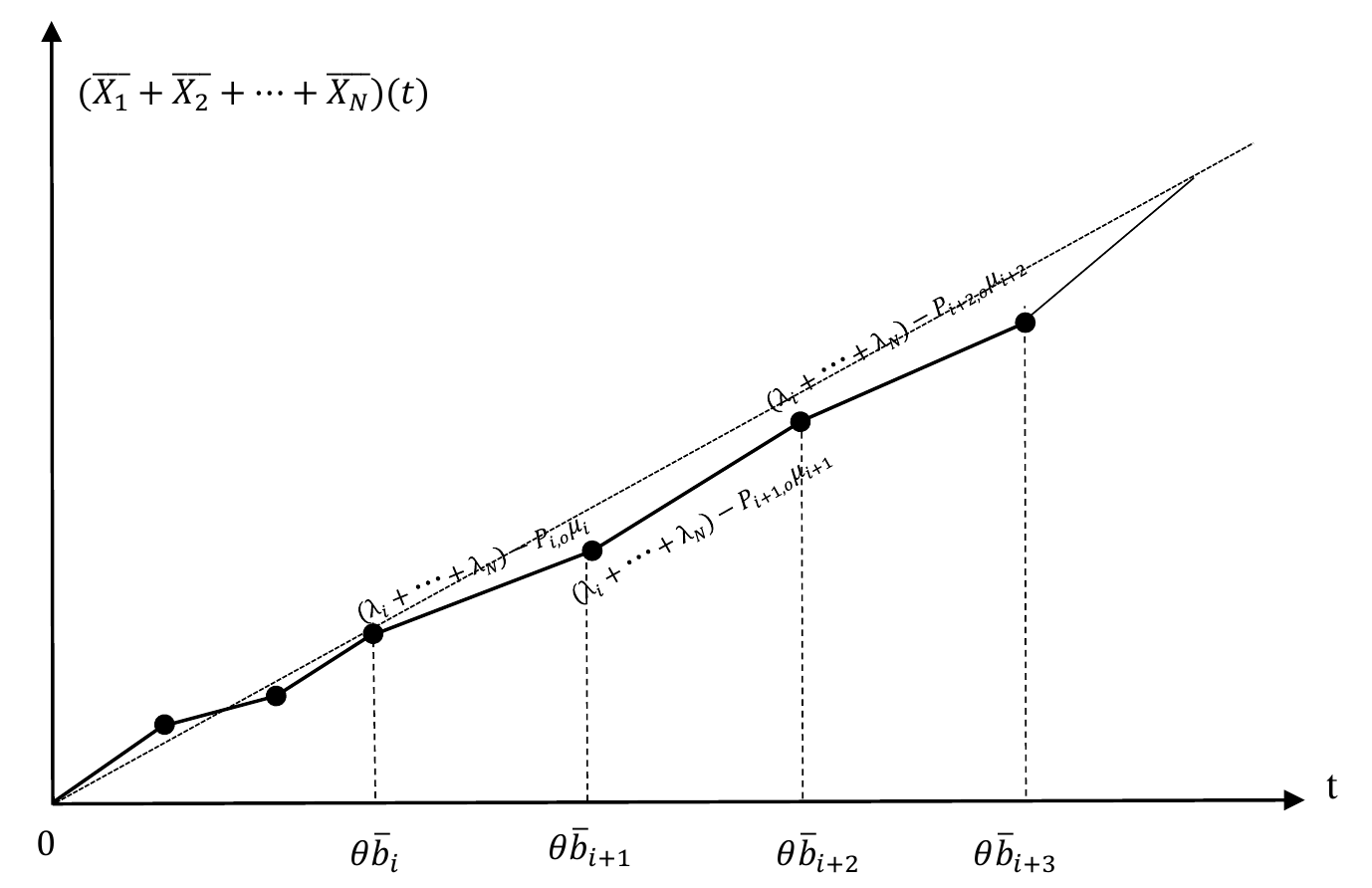}
\caption{\small Fluid limit of total population $\bar{X}_1 + \cdots + \bar{X}_N$. \;$\cdots \cdots$ is the average growth rate of total population.}
\label{fig:graph}
\end{figure}

According to Theorem \ref{thmcheckQ} and Corrollary \ref{cor_total}, the fluid limit processes both demonstrate an oscillation waveform with increasing amplitude and cycle time forward through time and oscillate at an infinite rate when approaching zero. To be more specific, the amplitude and cycle time both increase by $\theta-1$ times each cycle. Hence, the average growth rates of the fluid limit in each cycle equal to the average growth rate in the whole time. As shown in Fig.\ref{fig:graph}, the average growth rate, denoted by $\beta$, can be given by
\begin{equation*}
  \sum_{i=1}^N\frac{\sum_{j=1}^N\bar{a}_{i,j}+\sum_{j=1}^N\bar{a}_{i+1,j}}{2}(\bar{b}_{i+1}-\bar{b}_{i})=\int_{\bar{b}_{1}}^{\theta\bar{b}_{1}}\beta t dt,
\end{equation*}
which yields
\begin{equation}\label{averagerate}
  \beta=\frac{1}{\theta^2-1}\left[\sum_{i=1}^N\left(\sum_{j=1}^N\bar{a}_{i,j}+\sum_{j=1}^N\bar{a}_{i+1,j}\right)(\bar{b}_{i+1}-\bar{b}_{i})\right].
\end{equation}

By the definition of the scaled queue length process,  the fluid limit could approximate the original queue length process in steady state. Furthermore, the average growth rate in \eqref{averagerate} allows us to study the optimization problem of how to choose the gating indexes of each queue to minimize the total queue length.    Because each of the queues adheres to a branching-type service discipline, we also study  how to choose the exhaustiveness $f_i$ with the same objective in mind. Here we only provide an example by utilizing the genetic algorithm to solve the optimization problem in Section \ref{numerical}.

\section{Proof of Lemma \ref{exhaustiveness} and Lemma \ref{ineq-infty}}\label{prove-lemma}


\noindent{\bf Proof of Lemma \ref{exhaustiveness}.} By (4.1) and (4.2) in \cite{liu2015asymptotic}, we have
  $f_i=\frac{\E T_i}{\varphi_i}$, where $\varphi_i$ is defined in Lemma \ref{eq-variable}.
Hence, it only remains to prove \eqref{visit}.

For $k \in N \cup \{\infty\}$, let $T_i(k)$
be the visit duration at  $Q_i$ given that the service discipline at  $Q_i$ is $k$-gated.  Obviously, we have $\E T_i(\infty) = \frac{1}{\mu_i(1 - \frac{\lambda_i}{\mu_i}-p_{i,i})}$.
Since
\begin{eqnarray*}\label{eq-t2}
T_i(0) = 0 \quad \mbox{and} \quad T_i(k+1) \overset{d}{=} B_i + \sum_{l=1}^{E_i(B_i)} T_i^{(l)}(k) + \delta_{\{A_{i} = i\}} T_i(k), \qquad k \in \N,
\end{eqnarray*}
where $T_i^{(l)}(k), l \in \N^+$ are i.i.d. copies of $T_i(k)$ and the random elements $B_i$, $E_i(\cdot)$ and $\{T_i^{(l)}(k)\}_{l \in \N^+}$ are mutually independent, we get
\begin{eqnarray}\label{eq-exp-v}
\E T_i(k+1) & = & \frac{1}{\mu_i} + (\frac{\lambda_i}{\mu_i} + p_{i,i}) \E T_i(k)\nonumber\\
& = & \frac{1}{\mu_i} (1 + \frac{\lambda_i}{\mu_i} + p_{i,i}) + (\frac{\lambda_i}{\mu_i} + p_{i,i})^2 \E T_i(k - 1)\nonumber\\
& = & \cdots = \frac{1}{\mu_i} \frac{1-(\frac{\lambda_i}{\mu_i} + p_{i,i})^{k+1}}{1-\frac{\lambda_i}{\mu_i} - p_{i,i}}.
\end{eqnarray}
Therefore,
\begin{eqnarray*}
\E T_i = \sum_{k \in \N \cup \{\infty\}} \P(\kappa_i = k) \E T_i(k) =\frac{1 - \E (\frac{\lambda_i}{\mu_i} + p_{i,i})^{\kappa_i}}{\mu_i(1 - \frac{\lambda_i}{\mu_i}-p_{i,i})}.
\end{eqnarray*}

In order to prove Lemma \ref{ineq-infty}, we need the following lemmas.

\begin{lem}(\cite{remerova2013random}, Lemma 4)\label{def-eq-fx}
 Suppose that a function $f(\cdot): [0, \infty) \to [0, \infty)$ is bounded in a finite interval $[0,T]$ and nondecreasing in $[T, \infty)$, and that
 $f(x) \to \infty$ as $x \to \infty$. Suppose also that, for some (and, hence, for all) $c > 1$,
 \begin{eqnarray*}\label{de-fx}
 \limsup_{x \to \infty} \frac{f{(c x)}}{f(x)} < \infty.
 \end{eqnarray*}
 Consider an i.i.d. sequence $\{Y^{(n)}\}_{n \in \N^+}$ of nonnegative, non degenerate at $0$ random variables, and the renewal process
 \begin{eqnarray*}
 Y(t) = \max\{ n \in \N: \sum_{k=1}^{n} Y^{(k)} \leq t\}, \qquad t \in [0, \infty).
 \end{eqnarray*}
 Let $\tau$ be a nonnegative random variable which may depend on the sequence $\{Y_n\}_{n \in \N^+}$. Assume that $\E f(\tau) < \infty$. Then
 $\E f(Y(\tau))$ is finite too.
\end{lem}

\begin{lem}(\cite{remerova2013random}, Lemma 5)\label{eq-sum-fx}
Consider a sequence $Y^{(n)}_{n \in \N^+}$ of nonnegative random variables that are identically distributed (but not necessarily independent), and also a $\N$-valued random variable $\eta$ that does not depend on $Y^{(n)}_{n \in \N^+}$. If $f(\cdot): [0, \infty) \to \R$ is a convex function, then
\begin{eqnarray*}
\E f(\sum_{k=1}^{\eta}Y^{(k)}) \leq \E f(\eta Y^{(1)}).
\end{eqnarray*}
\end{lem}

We continue proving  Lemma \ref{ineq-infty}. It suffices to show that
\begin{eqnarray*}
\E f(L_{i j}) < \infty, \qquad \mbox{for all} \  i \ \mbox{and} \  j,
\end{eqnarray*}
where
\begin{eqnarray*}
f(x) = \left\{
         \begin{array}{ll}
           0, & \hbox{$x \in [0,1]$;} \\
           x \log x, & \hbox{$x \in [1, \infty)$.}
         \end{array}
       \right.
\end{eqnarray*}
Note that the function $f(\cdot)$ is convex in $(1, \infty)$, its derivative $\log(\cdot) + 1$ is nondecreasing, and in the other points, it is easy to verify the convexity. Note also that
\begin{eqnarray}\label{def-ineq-fx}
f(x y) \leq x f(y) + y f(x), \qquad x, y \in [0, \infty).
\end{eqnarray}

\noindent {\bf Proof of Lemma \ref{ineq-infty}.}
(1) {\bf Finiteness of $\E f(T_i)$.} Recall the definition of $T_i(k)$ in the proof of Lemma \ref{exhaustiveness}, we have
\begin{eqnarray*}
T_i(k)\uparrow T_i(\infty) \qquad a.s. \quad \mbox{as} \quad k \to \infty.
\end{eqnarray*}

Then by the monotonicity, convexity of $f(\cdot)$ and the auxiliary  Lemma \ref{eq-sum-fx} combined with \eqref{def-ineq-fx}, we have
\begin{align*}
\E f(T_i(k))  \leq &  \E f(T_i(k+1)) = \E f(B_i + \sum_{l=1}^{E_i(B_i)} T_i^{(l)}(k) + \delta_{\{A_{i} = i\}} T_i(k))\\
 \leq & \frac{1}{3} \E f(3 B_i) + \frac{1}{3}\E f(3\sum_{l=1}^{E_i(B_i)} T_i^{(l)}(k)) + \frac{1}{3} \E f(3 \delta_{\{A_{i} = i\}} T_i(k))\\
 \leq & \frac{1}{3} \E f(3 B_i) + \frac{1}{3}\E f(3 E_i(B_i) T_i^{(l)}(k)) + \frac{1}{3} \E f(3\delta_{\{A_{i} = i\}} T_i(k))\\
 \leq & \frac{1}{3} \E f(3 B_i) + \E  [E_i(B_i)] \E f(T_i^{(l)}(k)) + \frac{1}{3} \E T_i^{(l)}(k) \E f(3 E_i(B_i))\\
& + \E \delta_{\{A_{i} = i\}}\E f(T_i(k)) +　\frac{1}{3} \E T_i(k) \E f(3\delta_{\{A_{i} = i\}})\\
 \leq & \frac{1}{3} \E f(3 B_i) + (\frac{\lambda_i}{\mu_i}+p_{i,i}) \E f(T_i(k)) + \frac{1}{3} \E T_i(k)[f(3)+ \E f(3 E_i(B_i))].
\end{align*}
By Lemma \ref{def-eq-fx} and \eqref{eq-exp-v}, we have $\E f(3 E_i(B_i)) < \infty$ and $\E T_i(k) \leq \frac{1}{\mu_i(1 - \frac{\lambda_i}{\mu_i} - p_{i,i})}$.
Therefore, for all $k \geq 2$, we get
\begin{eqnarray*}
\E f(T_i(k)) \leq \frac{C}{1-\frac{\lambda_i}{\mu_i} -p_{i,i}},
\end{eqnarray*}
where
\begin{eqnarray*}
C= \frac{1}{3} \E f(3 B_i) + \frac{1}{3}\frac{f(3)+ \E f(3 E_i(B_i))}{\mu_i(1 - \frac{\lambda_i}{\mu_i} - p_{i,i})} < \infty.
\end{eqnarray*}
(2) { \bf Finiteness of $\E f(\check{L}_{i,i})$.} $\check{L}_{i,i}$ is bounded stochastically from above by the number of customers arriving from outside when the system is empty and during $T_i(\infty)$ (the visit duration at $Q_i$ with exhaustive service policy). Therefore, $\check{L}_{i,i} < 1 + E_i(T_i(\infty))$. Hence, $\E f(\check{L}_{i,i}) < \E f(1 + E_i(T_i(\infty))) < \infty$ by Lemma \ref{def-eq-fx}.

\noindent
(3) { \bf Finiteness of $\E f(L_{i,j})$.} The result can be proved in the same way as  in the proof of  Lemma 3 in \cite{remerova2013random}.

\section{Proof of Theorems \ref{thmQ} and  \ref{thmcheckQ}}\label{prove-theorem}
To prove Theorems \ref{thmQ} and  \ref{thmcheckQ}, we need some further notations:
\begin{enumerate}
  \item Define the renewal processes
$$Y_{i}(t) = \max\left\{n\in \mathbb{N}^+, \text{such that} \sum_{i=1}^n B_i^{(k)}\leq t\right\},$$
where $B_i^{(k)}$ are i.i.d. copies of $B_i$.
  \item Let $I_{i}(t)$ be the whole time that the server has spent at $Q_i$ before time $t$, i.e.,
    $$I_{i}(t) = \int_0^t I( \mbox{ queue}\; i\; \mbox{is in service in time } \;s ) ds \qquad t \in (0, \infty).$$
\item Define index $\nu$ by
\begin{equation*}
\nu= \max\left\{n\in \mathbb{N}^+, \text{such that} \ {\bf X}( t^{(n)}) = {\bf 0}\  \mbox{and}\  {\bf{X}}(t^{(m)}) \neq {\bf 0}\  \mbox{for all}\  m > n\right\}.
\end{equation*}
\end{enumerate}

\begin{lem}(\cite{remerova2013random}, Proposition 2)\label{eq-lim-yn}
Let a random variable $Y$ have a finite mean value and, for each $n \in \N^+$, let $Y_n^{(k)}, k \in \N^+$ be i.i.d. copies of $Y$. Let $\tau_n, n \in \N^+$ be $\N$-valued random variables such that $\tau_n$ is independent of the sequence $\{Y_n^{(k)}\}_{k \in \N^+}$ for each $n$ and $\tau_n \to \infty$ in probability as $n \to \infty$. Finally, let a sequence $\{T_n\}_{n \in \N^+}$ of positive numbers increase to $\infty$. If there exists an a.s. finite random variable $\tau$ such that $\tau_n/T_n \to \tau$ in probability as $ n \to \infty$, then
\begin{eqnarray*}
\sum_{k=1}^{\tau_n} \frac{Y_n^{(k)}}{T_n} \to \tau \E Y \qquad \mbox{in probability } \qquad \mbox{as} \quad n \to \infty.
\end{eqnarray*}
\end{lem}
\begin{lem}\label{lem-ti}
For $i = 1, \cdots, N$, there exist constants $b_i \in (0, \infty)$ and ${\bf a}_i = (a_{i,1}, \cdots, a_{i,N}) \in [0, \infty)^N$ such that
\begin{eqnarray*}
\frac{t_i^{(n)}}{\theta^n} \to b_i \xi \quad \mbox{and} \quad \frac{{\bf X} (t_i^{(n)})}{\theta^n} \to \xi {\bf a}_i \quad a.s. \quad \mbox{as } \quad n \to \infty.
\end{eqnarray*}
The $b_i$'s and ${\bf a}_i$'s are specified by
\begin{align}\label{eq-def-bi}
&b_1 = \frac{\sum_{i=1}^N v_i \bar{c}_i}{\rho- 1}, \nonumber\\
&b_{i+1} = b_i + \Big(v_i + \lambda_i(b_i - b_1) + \sum_{j=1}^{i-1} p_{j,i} \mu_j(b_{j+1} - b_j)\Big)t_i, \quad i=1,\cdots,N,
\end{align}
and
\begin{eqnarray} \label{eq-def-ai}
&&{\bf a}_1 = {\bf v},  a_{i+1, j} = \left\{
                                       \begin{array}{ll}
                                         a_{i,j} + (\lambda_j + p_{i, j} \mu_i)(b_{i+1} - b_{i}), & \hbox{$ j \neq i$;} \\
                                         a_{i,i} + (\lambda_i - \mu_i (1 - p_{i, i}))(b_{i+1} - b_i), & \hbox{$ j = i$.}
                                       \end{array}
                                     \right.
\end{eqnarray}
Obviously, the ${\bf a}_i$'s also satisfy
\begin{eqnarray}\label{eq-def-ai-also}
&&{\bf a}_1 = {\bf v}, \qquad {\bf a}_{i+1} = {\bf a}_{i} - a_{i,i} {\bf e}_i+ a_{i,i} {\bf \check{m}}_{i}={\bf a}_{i}{\bf M}_i \quad i = 1, \cdots, N.
\end{eqnarray}
\end{lem}

{\bf Proof } (1) {\bf Limit of $t_1^{(n)}/\theta^n$.} By the definition of $\nu$, which is a.s. finite, combined with the total workload process, we have, for $n > \nu$,
\begin{eqnarray*}\label{eq-def-t1}
t_1^{(n)} = t^{(n)} = \sum + \sum_{i=1}^{N}\sum_{k=1}^{E_i(t^{(n)})}\bar{B}_i^{(k)} - \sum_{i=1}^{N}\sum_{k=1}^{X_i(t^{(n)})} \bar{B}_i^{(k)} ,
\end{eqnarray*}
where $\bar{B}_i^{(k)}$ are i.i.d. copies of $\bar{B}_i$ and $\sum = \sum_{l=0}^\nu (t_1^{(l)} - t^{(l)})$ is a.s. finite.  Therefore, we have
\begin{eqnarray*}\label{def-eq-t1n}
t_1^{(n)} = t^{(n)} = \sum + t^{(n)}A_1^{(n)} -  \theta^n A_{2}^{(n)},
\end{eqnarray*}
where
\begin{eqnarray*}
&& A_1^{(n)} = \sum_{i=1}^{N} \frac{\sum_{k=1}^{E_i(t^{(n)})} \bar{B}_i^{ (k)}}{E_i(t^{(n)})}\frac{E_i(t^{(n)})}{t^{(n)}},\\
&& A_2^{(n)} = \sum_{i=1}^{N} \frac{\sum_{k=1}^{X_i(t^{(n)})} \bar{B}_i^{(k)}}{X_i (t^{(n)})}\frac{X_i(t^{(n)})}{\theta^n},
\end{eqnarray*}
then
\begin{eqnarray}\label{eq-tran-t1}
\frac{t_1^{(n)}}{\theta^n} = \frac{t^{(n)}}{\theta^n} = \frac{A_2^{(n)} - \sum/\theta^n}{A_1^{(n)} - 1}.
\end{eqnarray}
By the SLLN and Proposition \ref{pro-limit}, we obtain, as $n \to \infty$,
\begin{align}\label{eq-lim-A12}
A_1^{(n)} \to \sum_{i=1}^{N} \lambda_i \E \bar{B}_i^{(k)} = \sum_{i=1}^{N}\lambda_i \bar{c}_i=\rho,\quad
 A_2^{(n)} \to \sum_{i=1}^{N} v_i \xi \E \bar{B}_i^{ (k)} = \sum_{i=1}^{N} v_i \bar{c}_i \xi\; \; a.s..
\end{align}
Then by \eqref{eq-tran-t1} and \eqref{eq-lim-A12}, we have, as $n \to \infty$,
\begin{eqnarray*}\label{eq-lim-tn-b1}
\frac{t^{(n)}}{\theta^n} \to b_1 \xi \quad \mbox{and} \quad \frac{t_1^{(n)}}{\theta^n} \to b_1 \xi,
\end{eqnarray*}
where $b_1 = \frac{\sum_{i=1}^{N} v_i \bar{c}_i}{\rho - 1}$.

(2) {\bf  Limit of $t_i^{(n)}/\theta^n$.} In (1), by utilizing the index $\nu$ and equation $t_1^{(n)}=t^{(n)}$, we proved $\lim_{n \to \infty}t_1^{(n)}/\theta^n = b_1 \xi$. By the symmetry, there should also exist positive numbers $b_i$ such that
\begin{eqnarray*}\label{eq-lim-tin}
\lim_{n \to \infty} \frac{t_i^{(n)}}{\theta^n}= b_i \xi, \; \; i=1, \cdots, N.
\end{eqnarray*}
It remains to prove these positive numbers $b_i$'s satisfy \eqref{eq-def-bi}, which refer to (5) below.

(3){ \bf Limit of $X(t_i)/ \theta^n$ and \eqref{eq-def-ai}. } By definition, we have
{\small
\begin{equation}\label{eq-def-QJ}
X_j(t_{i+1}^{(n)}) =\left\{
                      \begin{array}{ll}
                         X_j(t_i^{(n)}) + E_j(t_{i+1}^{(n)}) - E_j(t_i^{(n)}) + \sum\limits _{k=1}^{Y_i(I_i(t^{n+1})) - Y_i(I_i(t^{n}))} \delta_{\{A_{i}^{(k)} = j\}}, & \hbox{$j \neq i$;} \\
                         X_i(t_i^{(n)}) + E_i(t_{i+1}^{n}) - E_i(t_{i}^{(n)})- \sum\limits_{k=1}^{Y_i(I_i(t^{n+1})) - Y_i(I_i(t^{n}))} \delta_{\{A_{i}^{(k)} \neq i\}}, & \hbox{$j = i$,}
                      \end{array}
                    \right.
\end{equation}
}
where $A_{i}^{(k)}$ are i.i.d. copies of $A_{i}$.
Therefore, by SLLN and \eqref{eq-lim-tin}, we obtain
\begin{eqnarray*}\label{eq-lim-Ej}
\frac{E_j(t_{i+1}^{(n)}) - E_j(t_i^{(n)})}{\theta^n} \to \lambda_j(b_{i+1} - b_i) \xi \qquad a.s. \quad \mbox{as} \quad n \to \infty,
\end{eqnarray*}
and
\begin{align*}\label{eq-lim-I}
&\frac{\sum_{k=1}^{Y_i(I_i(t^{n+1})) - Y_i(I_i(t^{n}))} \delta_{\{A_{i}^{(k)} = j\}}}{\theta^n}\nonumber\\
&=  \frac{\sum_{k=1}^{Y_i(I_i(t^{n+1})) - Y_i(I_i(t^{n}))} \delta_{\{A_{i}^{(k)} = j\}}}{Y_i(I_i(t^{n+1})) - Y_i(I_i(t^{n}))} \frac{Y_i(I_i(t^{n+1})) - Y_i(I_i(t^{n}))}{ I_i(t^{(n + 1)}) - I_i(t^{(n)})} \frac{I_i(t^{n+1}) - I_i(t^{n})}{\theta^n}\nonumber\\
&\to \mu_i (b_{i+1} - b_{i}) \xi \E \delta_{\{ A_{i}^{(k)} = j\}}  = p_{i,j} \mu_i(b_{i+1} - b_i) \xi.
\end{align*}
Similarly,
\begin{equation*}
  \frac{\sum_{k=1}^{Y_i(I_i(t^{n+1})) - Y_i(I_i(t^{n}))} \delta_{\{A_{i}^{(k)} \neq i\}}}{\theta^n}
\to (1-p_{i,i}) \mu_i(b_{i+1} - b_i) \xi.
\end{equation*}
It follows that
\begin{eqnarray*}
\frac{X_j(t_{i+1}^{(n)})}{\theta^n} \to a_{i+1, j}\xi,
\end{eqnarray*}
where
\begin{equation*}\label{eq-def-aij}
a_{i+1, j} =\left\{
              \begin{array}{ll}
                a_{i,j} + \lambda_j (b_{i+1} - b_i) + p_{i,j} \mu_j (b_{i+1} - b_i), & \hbox{$j \neq i$;} \\
                a_{i,i} + \lambda_i(b_{i+1} - b_{i}) - \mu_i(1-p_{i,i})(b_{i+1} - b_i), & \hbox{$j = i$.}
              \end{array}
            \right.
\end{equation*}

(4) {\bf Equation \eqref{eq-def-ai-also}.} \eqref{eq-def-ai-also} follows from the below equation
\begin{eqnarray*}
{\bf Q}(t_{i+1}^{(n)}) = {\bf Q}(t_i^{(n)}) - Q_i(t_i^{(n)}){\bf e}_i + \sum_{k=1}^{Q_i(t_i^{(n)})}{ \bf \check{L}}_i^{(n,k)},
\end{eqnarray*}
where ${\bf \check{L}}_{i}^{(n,k)} = (\check{L}_{i,1}^{(n,k)}, \cdots, \check{L}_{i,N}^{(n,k)})$ and $\check{L}_{i,j}^{(n,k)}$ are i.i.d. copies of $ \check{L}_{i,j} $.

(5) {\bf Equation \eqref{eq-def-bi}.}
Recall that
\begin{align}
&t_{i+1}^{(n)} = t_i^{(n)} + \sum_{k=1}^{X_i(t_i^{(n)})} T_i^{(k)}, \label{eq-def-tin}\\
&X_i(t_i^{(n)}) = X_i(t_1^{(n)}) + E_i(t_i^{(n)}) - E_i(t_1^{(n)}) + \sum_{j=1}^{i-1} \sum_{k=1}^{Y_j (I_j(t^{(n)})) - Y_j (I_j(t^{(n-1)}))} \delta_{\{A_{j}^{(k)} = i\}}\label{eq-def-qitinn}.
\end{align}
From Lemma \ref{eq-lim-yn} and \eqref{eq-def-tin}, we get
\begin{eqnarray}\label{eq-def-lim-ti}
b_{i+1} - b_i = a_{i,i} t_i.
\end{eqnarray}
By \eqref{eq-def-qitinn} and the SLLN, we obtain
\begin{eqnarray}\label{eq-def-lim-aii}
a_{i,i} = v_i + \lambda_i(b_i - b_1) + \sum_{j=1}^{i-1}  \mu_j(b_{j+1} - b_{j}) p_{j,i}.
\end{eqnarray}
Then \eqref{eq-def-bi} can be proved by substituting \eqref{eq-def-lim-aii} into \eqref{eq-def-lim-ti}.

By $t_{N+1}^{(n)} = t_1^{(n + 1)}$, we obtain $b_{N+1} \xi = \lim\limits_{n \to \infty} t_{N+1}^{(n)}/\theta^n = \lim\limits_{n \to \infty}\theta  t_1^{(n+1)}/\theta^{n+1} = \theta b_1 \xi$, i.e. $b_{N+1}= \theta b_1$. This can be proved as follows. By the definition of ${\bf M}_i$ in Lemma \ref{eq-variable}, it is easy to give
\begin{equation*}\label{right}
  ({\bf M}_i-{\bf I})\bar{\bf c}^T=(\rho-1)t_i {\bf e}_i,
\end{equation*}
where ${\bf I}$ is the identity matrix and $\bar{\bf c}=(\bar{c}_1,\ldots,\bar{c}_N)$.
Substituting the above equation into \eqref{eq-def-lim-ti} yields
\begin{equation*}\label{right}
  b_{i+1} - b_i = a_{i,i} t_i={\bf a}_i{\bf e}_it_i=\frac{1}{\rho-1}{\bf a}_i({\bf M}_i-{\bf I})\bar{\bf c}^T=\frac{1}{\rho-1}\left({\bf a}_{i+1}-{\bf a}_i\right)\bar{\bf c}^T.
\end{equation*}
which gives immediately
\begin{align*}\label{right}
  b_{N+1}=&\sum_{i=1}^N(b_{i+1} - b_i)+b_1 =\frac{1}{\rho-1}\left({\bf a}_{N+1}-{\bf a}_1\right)\bar{\bf c}^T+b_1\\
  =&\frac{1}{\rho-1}\left({\bf v}{\bf M}-{\bf v}\right)\bar{\bf c}^T+b_1=(\theta-1)\frac{{\bf v}\bar{\bf c}^T}{\rho-1}+b_1=(\theta-1)b_1+b_1=\theta b_1.
\end{align*}

{\bf Proof of Theorem \ref{thmQ}.}
With the results of Lemma \ref{lem-ti}, Theorem \ref{thmQ} immediately follows from Lemma 6 in \cite{remerova2013random}.

{\bf Proof of Theorem \ref{thmcheckQ}.}
For each $i$, by \eqref{def-eq-barQ}, we know that the function $\bar{X}_i(\cdot)$ might have discontinuities only at $t = 0$ and $t = \theta^{k} \bar{b}_{i}$ for each $k\in \N$. Since the function $\bar{X}_i(\cdot)$ is c\`{a}dl\`{a}g, the continuity of $\bar{X}(\cdot)$ is evident in combination with the definition of $\bf{a}_i$.

Additionally, the uniform convergence on compact sets can be proved in the same way as in the proof of Theorem 2 in \cite{remerova2013random}. Hence, it suffices to prove the point-wise convergence \eqref{def-eq-barQ} for each $i=1,2,\ldots,N$.

By Theorem \ref{thmQ}, we have, as $ n \to \infty$,
\begin{eqnarray*}
&&\frac{t_i^{(\eta_n + k)}}{\theta^n} \to \theta^k \bar{b}_i \xi \quad \mbox{and} \quad \frac{{ \bf X}(t_i^{(\eta_n + k)})}{\theta^n} \to \xi \theta^k { \bf \bar{a}}_i,\\
&&\frac{I_i(t_i^{(\eta_n + k)})}{\theta^n} \to \theta^k \frac{\bar{b}_{i+1} - \bar{b}_i}{\theta(\theta -1)}\xi,\\
&&\frac{E_i(t)}{t} \to \lambda_i \quad \mbox{and} \quad \frac{Y_{i} (t)}{t} \to \mu_i.
\end{eqnarray*}
 We will show that,  as $n \to \infty$,
\begin{eqnarray}\label{eq-def-lim-barq}
{\bf \bar{X}}^{(n)}(t) \to \xi {\bf \bar{X}}(\frac{t}{\xi}) \quad \mbox{for all } \quad t \in [0, \infty),
\end{eqnarray}
 where ${\bf \bar{X}}^{(n)}(\cdot)$ is given by the form of \eqref{def-sal-leng}.

For $t=0$, the convergence of \eqref{eq-def-lim-barq} holds since the system starts empty. For each $i=1,\ldots,N$, if $t\in[\theta^k \bar{b}_i, \theta^k \bar{b}_{i+1})$, it remains to prove
 \begin{eqnarray*}
\bar{X}_j(t) =\left\{
                \begin{array}{ll}
                  \theta^k \bar{a}_{i,i} + [\lambda_i - \mu_i(1-p_{i,i})](t - \theta^k \bar{b}_i), & \hbox{$j=i$;} \\
                  \theta^k \bar{a}_{i,j} + [\lambda_j + \mu_i p_{i,j}](t - \theta^k \bar{b}_i), & \hbox{ $j\neq i$.}
                \end{array}
              \right.
\end{eqnarray*}

For all $n$ big enough, $\frac{t_i^{(\eta_n + k)}}{\theta^n} < t < \frac{t_{i+1}^{(\eta_n + k)}}{\theta^n}$
implying that $Q_i$ is in service during $[t_i^{(\eta_n + k)}, \theta^n t)$.
Hence,
{\footnotesize
\begin{equation*}
X_j(\theta^n t) =\left\{
                      \begin{array}{ll}
X_i(t_i^{(\eta_n+k)}) + E_i(\theta^n t) - E_i(t_i^{(\eta_n+k)}) - \sum\limits_{k=1}^{Y_i(I_i(\theta^n t)) - Y_i(I_i(t^{(\eta_n + k)}))} \delta_{\{A_{i}^{(k)} \neq i\}}, &\hbox{$j = i$;}\\
                         X_j(t_i^{(\eta_n + k)})+ E_j(\theta^n t)-E_j(t_i^{(\eta_n + k)})
+\sum\limits_{l=1} ^{Y_i(I_i(\theta^n t)) - Y_i(I_i(t^{(\eta_n + k)}))}\delta_{\{A_{i}^{(l)}= j\}}, &\hbox{$j \neq i$.}
                      \end{array}
                    \right.
\end{equation*}
}
Therefore,
\begin{equation*}
\bar{X}_j^{(n)}(t)=\frac{X_j(\theta^n t)}{\theta^n} \to \left\{
                      \begin{array}{ll}
\xi \theta^k \bar{a}_{i,i} + [\lambda_i - \mu_i(1- p_{i,i})](t -\theta^k \bar{b}_i \xi), & \hbox{$j = i$;}\\
                      \xi \theta^{k}\bar{a}_{i,j}+(\lambda_j+\mu_i p_{i,j})(t-\theta^k \bar{b}_{i} \xi), & \hbox{$j \neq i$,}
                      \end{array}
                    \right.
\end{equation*}
where the right hand-side actually equals $\xi \bar{X}_j(\frac{t}{\xi})$.
And the proof is concluded.

\section{Numerical validation and optimization of gating indexes}\label{numerical}
\subsection{Numerical validation}
\begin{table}[htbp]
\centering
 \caption{\label{tab.1}Parameter values in 3-queue polling network}
 \begin{tabular}{ll}
  \toprule
  Parameter& Considered parameter values \\
  \midrule
  Arriving rate & $\lambda_1=\lambda_2=\lambda_3=1$ \\
  Service rate & $\mu_1=8, \mu_2=5, \mu_3=2$\\
   &\multirow{3}{*}{ $P=(p_{i,j})_{3\times 3}=\left(
                                \begin{array}{ccc}
                                  0.1 & 0.25 & 0.2 \\
                                  0.2 & 0.1 & 0.2 \\
                                  0.2 & 0.1 & 0.25 \\
                                \end{array}
                              \right)$}\\
   Transition probability &\\
   &\\
  \bottomrule
 \end{tabular}
\end{table}

 This subsection is devoted to test  the validity of the fluid limits of the scaled queue length processes defined in \eqref{def-sal-leng}. For simplicity, we consider a 3-queue polling system described in Tab.\ref{tab.1} with exponentially distributed service times. For this model, it is readily to obtain $\rho_1=0.4749,~\rho_2=0.5194,~\rho_3=0.8625$ and $\rho=1.8568$, which belongs to the overloaded traffic case studied in this paper.

\begin{figure}[htp]
\centering
\includegraphics[width=2.6 in]{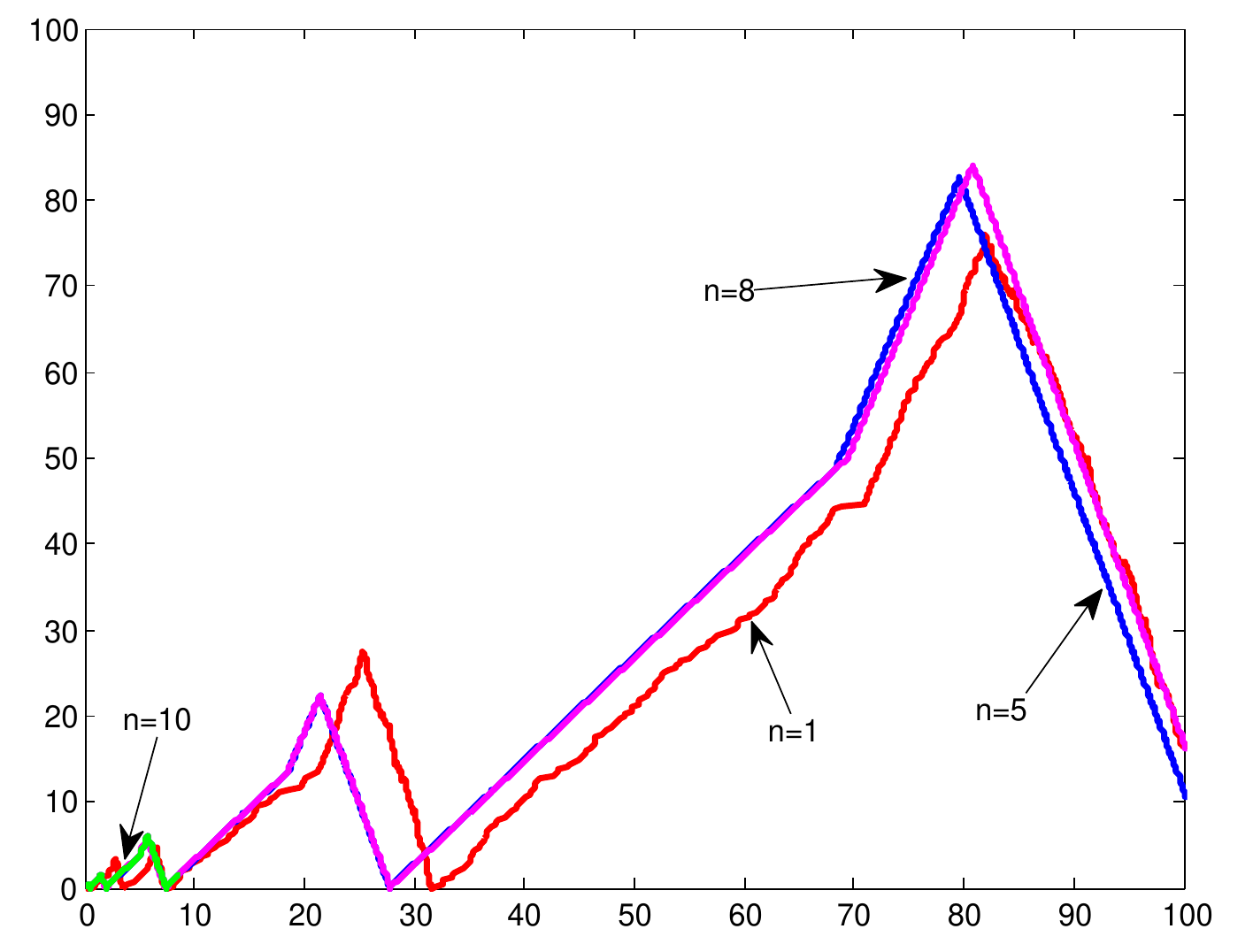}
\includegraphics[width=2.6 in]{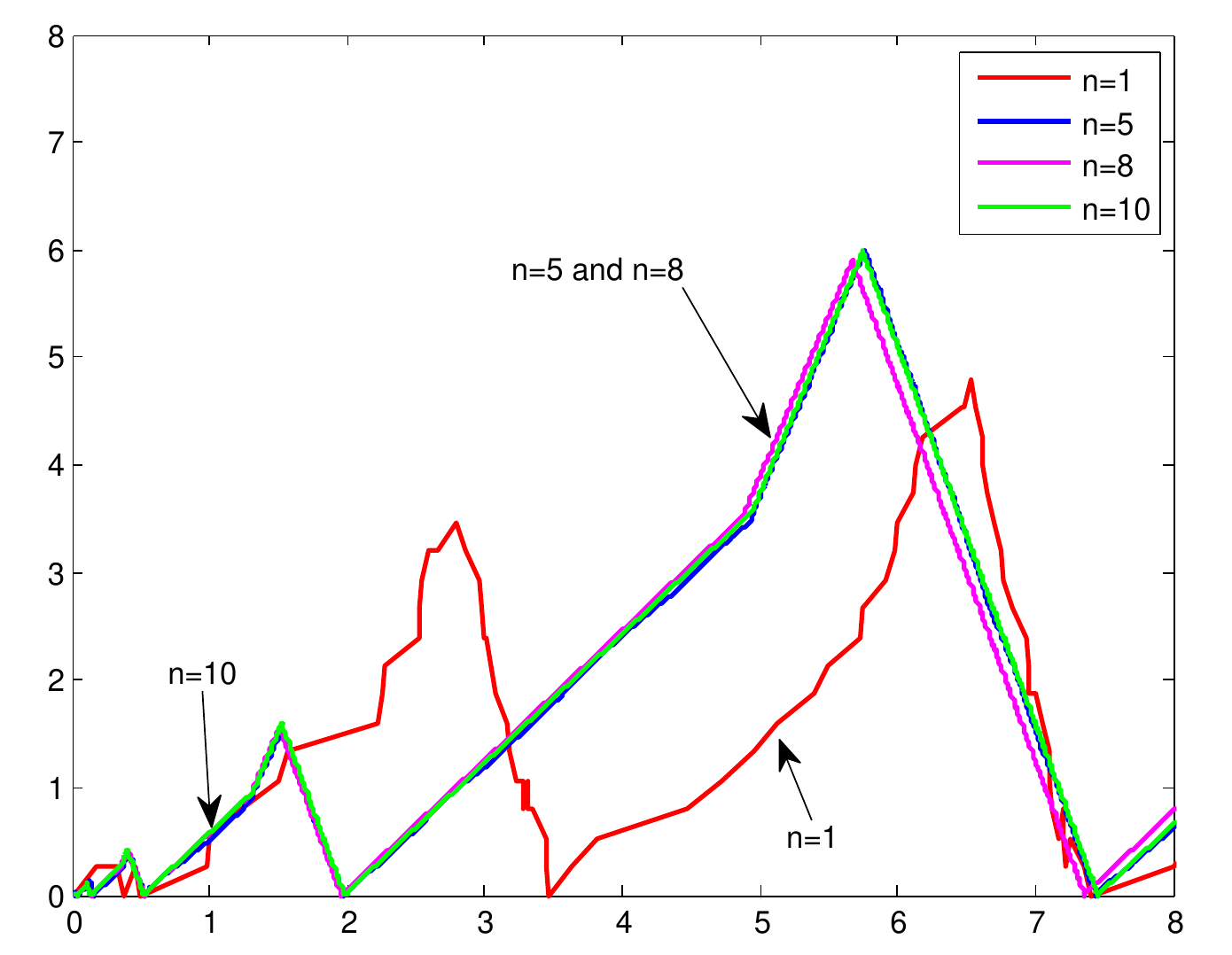} 
\caption {\small The scaled queue length process at $Q_2$ for different $n$ (left:$t\in [0,100]$, right:$t\in [0,8]$)}
\label{fig.3}
\end{figure}
\begin{figure}[htp]
\centering
\includegraphics[width=2.6 in]{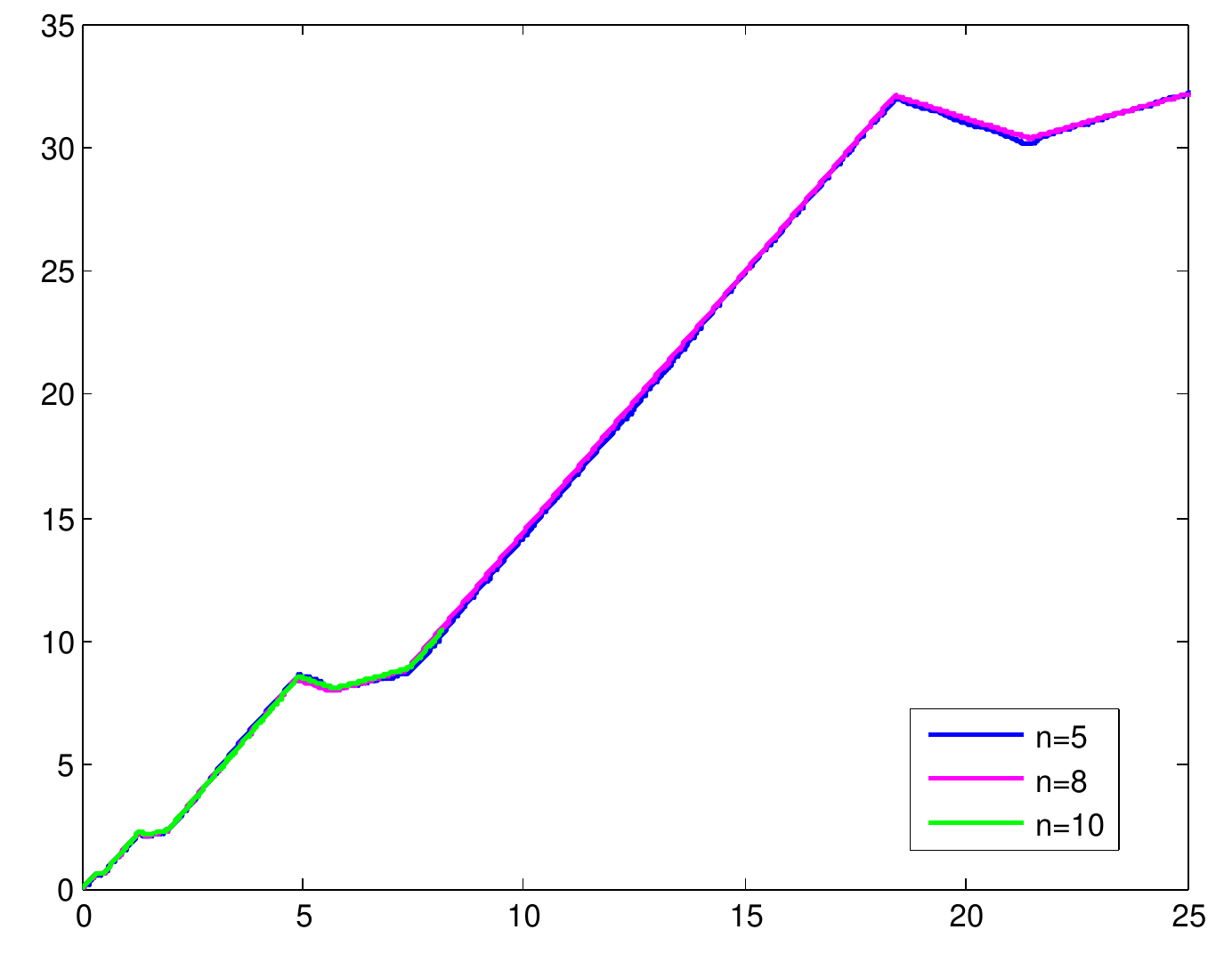}
\includegraphics[width=2.6 in]{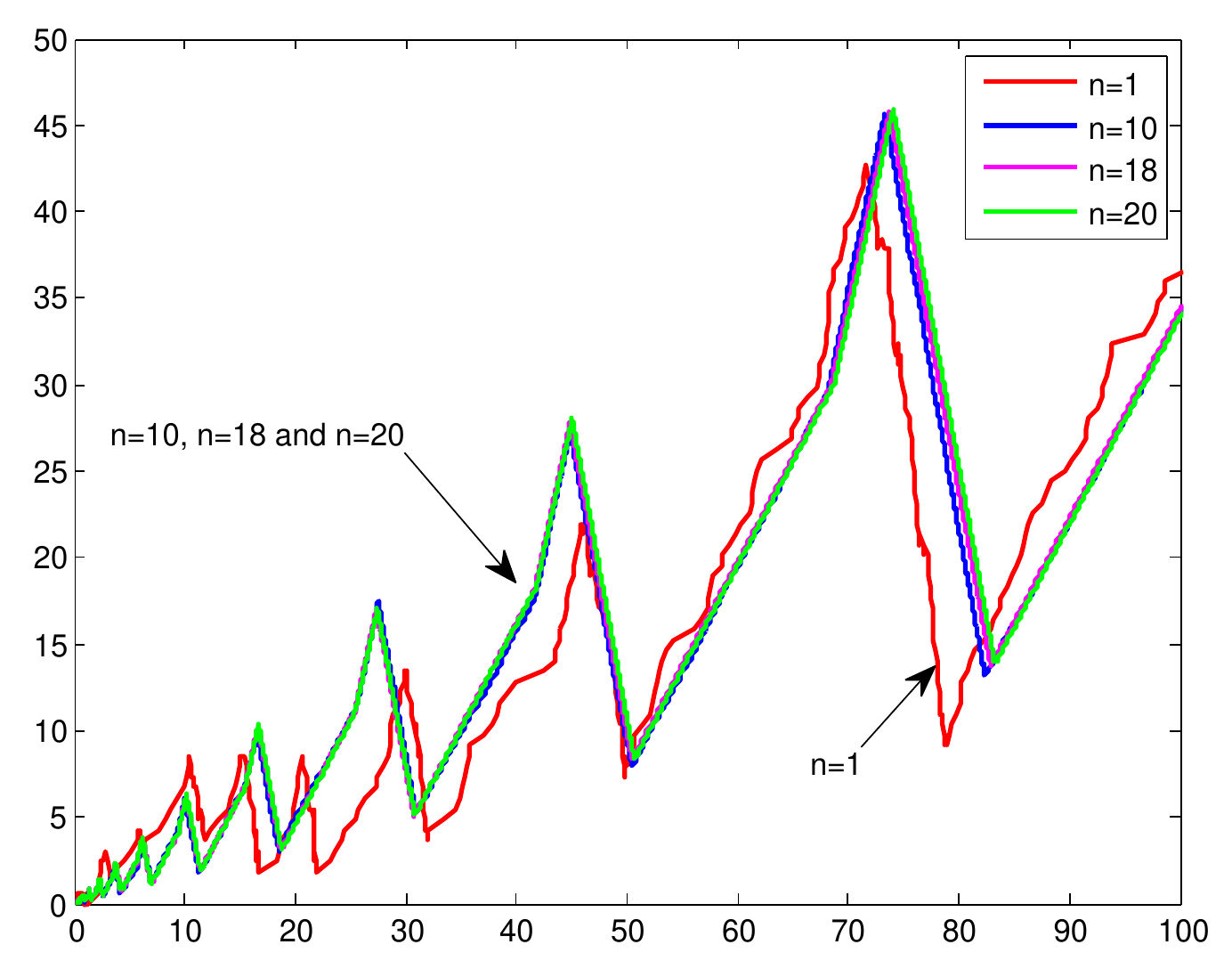} 
\caption {\small Left: the scaled total queue length process for different $n$ with exhaustive service policy\\
Right: the scaled queue length process at $Q_2$ for different $n$ with gated service policy}
\label{fig.4}
\end{figure}

 We utilize the SimEvents toolbox of Matlab to undertake the simulations of the polling networks. The exhaustive and gated service policies are taken for example and some vital variables are given in Tab.\ref{tab.2}. In order to illustrate the convergence of the scaled queue length processes defined in \eqref{def-sal-leng}, we take $n=1,5,8,10$ in polling network with exhaustive service policies and $n=1,10,18,20$ in the gated counterpart. The corresponding scaled queue length processes at $Q_2$ and the scaled total queue length process are depicted in Fig.\ref{fig.3} and Fig.\ref{fig.4}, respectively. Apparently, the scaled queue length sample paths get closer and closer as $n$ increases.

\begin{table}[htbp]
\centering
 \caption{\label{tab.2}Essential variable values in 3-queue polling network}
 \begin{tabular}{lll}
  \toprule
  \multirow{2}{*}{Variable}& \multicolumn{2}{c}{values} \\
  \cline{2-3}
  & Exhaustive & Gate\\
  \midrule
  Gating index & $\infty$ & $1$\\
  Exhaustiveness & $1$ & $1-(\frac{\lambda_i}{\mu_i} + p_{i,i})$\\
  Maximum eigenvalue & $\theta=3.7497$ & $\theta=1.6394$\\
  Left eigenvector &$v=[0.9731,0.683,0]$ & $v=[0.7454,0.5301,0.4774]$\\
  \bottomrule
 \end{tabular}
\end{table}

Moreover, as shown in Fig.\ref{fig.3} and Fig.\ref{fig.4}, the fluid limit processes both demonstrate an oscillation waveform with increasing amplitude and cycle time forward and oscillate at an infinite rate when approaching zero. According to Theorem \ref{thmcheckQ}, the amplitude and cycle time increase by $\theta-1$ times each cycle, which has been easily verified by the sample paths.

\subsection{Optimization of gating indexes}
Subsequently, we consider the optimization of the gating indexes by numerical method. It is assumed that the gating indexes are integers additionally. Virtually, the fluid limits only depend on the exhaustiveness of the service discipline at each queue (moment of gating index), which allows us to minimize the total queue length process through the accommodation of the integer gating indexes.

\begin{figure}[htp]
\centering
\includegraphics[width=2.6 in]{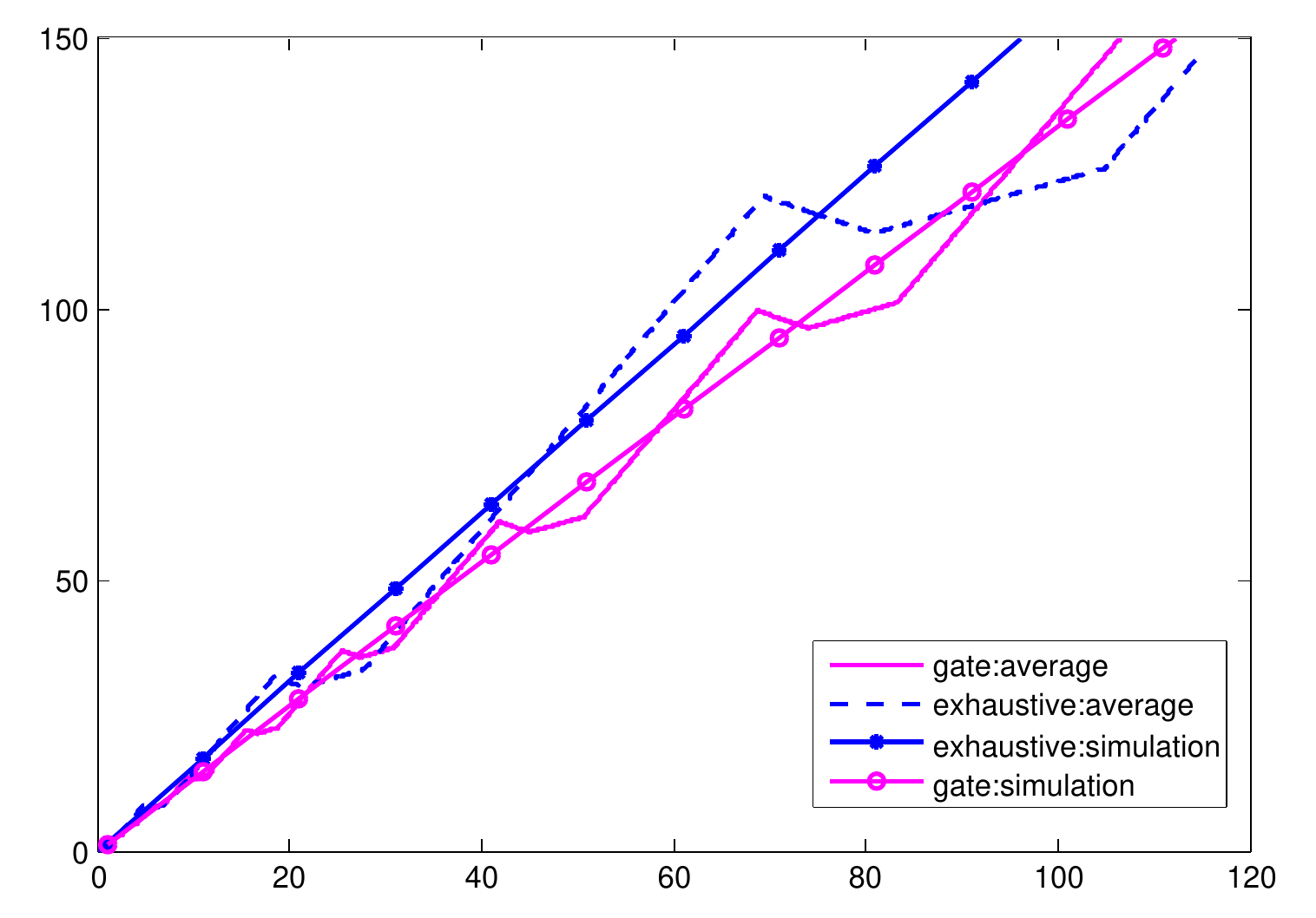}
\includegraphics[width=2.6 in]{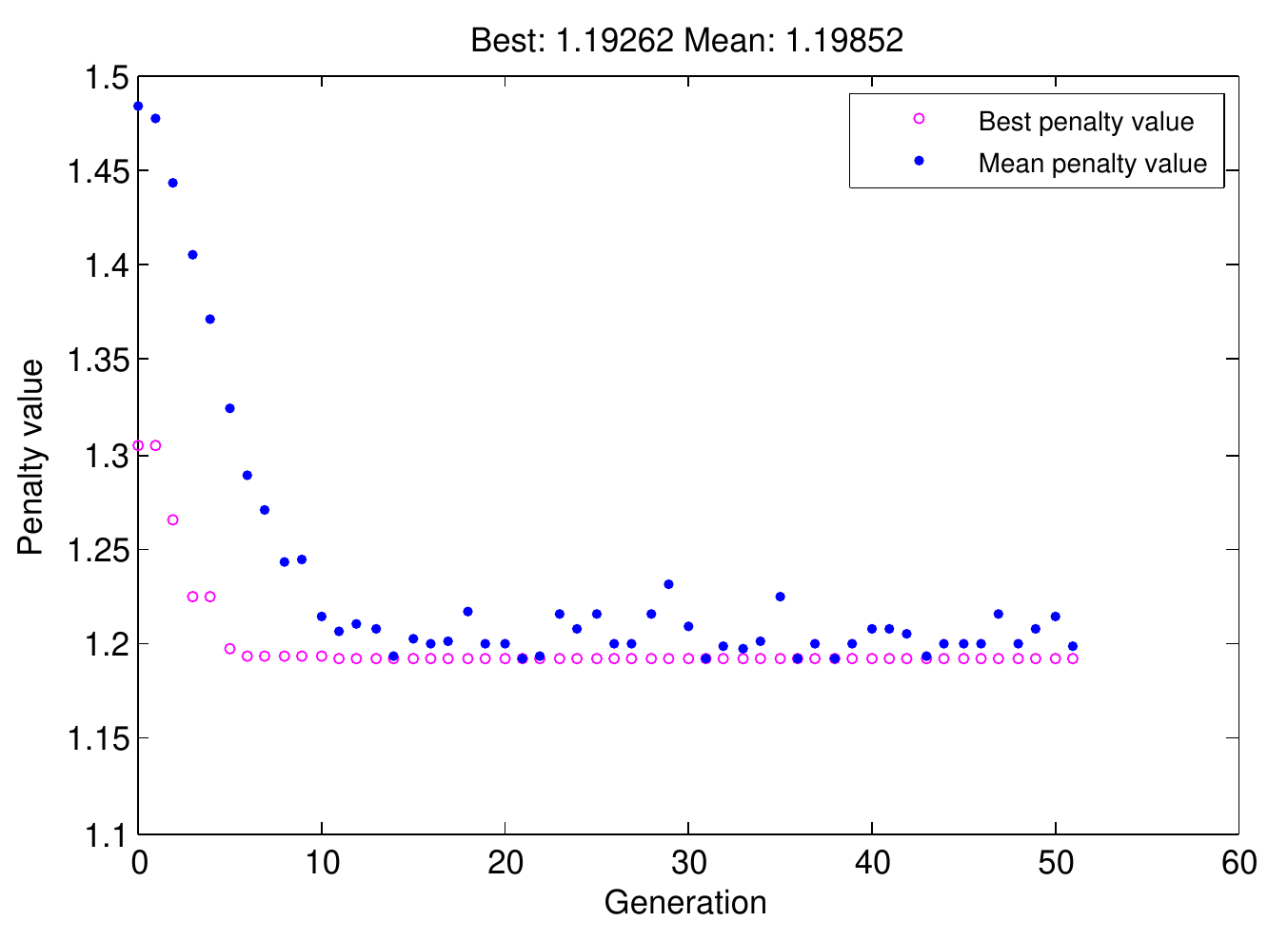} 
\caption {\small Left: the fluid limit of the scaled total queue length process (exhaustive and gate)\\
Right: the convergence process of the optimal average growth rate by GA solver}
\label{fig.5}
\end{figure}

By \eqref{averagerate}, the average growth rate of the total queue length process $\beta$ with exhaustive and gated service policies equals to 1.5025 and 1.2416 respectively (see Fig.\ref{fig.5}). This can be intuitively interpreted from the growth rate during each visiting period on different queues. The visiting period at $Q_3$ with the maximal growth rate (minimal service rate) takes 4 times as much time as others in exhaustive service policy. Instead, it takes less than 2 times as much time as others in gated service policy. Therefore, to minimize the average growth rate, we need to increase the visiting time at $Q_1$ and $Q_2$ and decrease the visiting time at $Q_3$.

To optimize the gating indexes turns to be an integer programming with three variables. The GA toolbox of Matlab is undertaken here to search for the optimal gating indexes. For our model, the GA solver just takes 51 iterations to find the optimal solution: $Q_1$ and $Q_2$ both take exhaustive service policy while $Q_3$ takes gated service policy. The minimal average growth rate equals to 1.19262 and the corresponding exhaustiveness is $f_1=f_2=1$, $f_3=0.25$. Fig.\ref{fig.5} depicts the process of the optimal average growth rate in each generation. Apparently, the convergence process turns to be very effective. Hence, the average growth rate of the fluid limit provides a simple and transparent method to optimize the gating index.

\section{Conclusions and Further Research}\label{conclude}

Inspired by \cite{remerova2013random}, we present the fluid limit of an overloaded polling system with general random multi-gated service discipline and customer re-routing policies. These results provide new fundamental insight in the impact of exhaustiveness. As a by-product, we propose an optimization problem of gating indexes to minimize the total queue length process.

This work gives rise to a variety of directions for further research. A logical follow-up step would be to study the case with non-zero switch-over time and more general branching-type polling models. In addition, the asysmptotic behaviors of discrete-time polling systems are also direct extensions to this study. Furthermore, the fluid limit allows us to propose control strategies of the growth depression, which requires substantially more effort.

\section*{Acknowledgement}
This research is partially supported by the National Natural Science Foundation of China (11671404, 11571052), by the Fundamental Research Funds for the Central Universities of Central South University(2015zzts012), and by the Fundamental Research Funds for the Central Universities(WUT:2017 IVA 069). 


\end{document}